\documentclass[10pt]{amsart}
\usepackage{amssymb}
\usepackage{amsmath}
\newtheorem{thm}{Theorem}[section]
\newtheorem{prop}[thm]{Proposition}
\newtheorem{cor}[thm]{Corollary}
\newtheorem{lem}[thm]{Lemma}
\newtheorem{rem}[thm]{Remark}

\newcommand{\C}{{\mathcal C}}

\newcommand\rep{\operatorname{rep}}
\newcommand\Irr{\operatorname{Irr}}
\newcommand\FPdim{\operatorname{FPdim}}

\newcommand\Hom{\operatorname{Hom}}

\newcommand\Pic{\operatorname{G}}
\newcommand{\1}{\textbf{1}}

\begin{document}

\title[On integral fusion categories]{On integral fusion categories with low-dimensional simple objects}
\author[J. Dong]{Jingcheng Dong}
\author[L. Dai]{Li Dai}

\address{College of Engineering, Nanjing Agricultural University, Nanjing 210031, Jiangsu, China}



\date{\today}

\begin{abstract}
Let $k$ be an algebraically closed
field of characteristic zero. In this paper we consider an integral fusion category over $k$ in which the Frobenius-Perron dimensions of its simple
objects are at most $3$.  We prove that such fusion category is of Frobenius type. In addition, we also prove that such fusion category is not simple.
\end{abstract}

\thanks{The research was partially
supported by the NSF of China (11201231), the China postdoctoral science foundation (2012M511643), the Jiangsu planned projects
for postdoctoral research funds (1102041C) and the agricultural
machinery bureau foundation of Jiangsu province (GXZ11003).}

 \maketitle

\section{Introduction}
The present work was motivated by an observation that the Frobenius-Perron dimensions of simple objects of many examples are small \cite{DNV,Nat2}. It was also motivated by a list of questions posed by a conference titled ``Classifying fusion categories" \cite{AMI}.

Let $k$ be an algebraically closed field of characteristic zero. A fusion category over $k$ is a $k$-linear semisimple rigid tensor category $\C$ which has finitely many isomorphism classes of simple objects, finite-dimensional hom spaces, and the unit object $\1$ of $\C$ is simple.

Fusion categories form a large class of categories. For example, if $G$ is a finite group, then the category $\rep G$ of its finite-dimensional representations is a fusion category over $k$. More generally, if $H$ is a finite-dimensional semisimple Hopf algebra over $k$, then the category $\rep H$ of its finite-dimensional representations is a fusion category. We refer the reader to \cite{ENO} for the main notions about fusion categories used throughout.

A fusion category $\C$ is called of Frobenius type if every Frobenius-Perron
dimension of its simple object divides the Frobenius-Perron
dimension of $\C$. Namely, the ratio $\FPdim
\C/\FPdim X$ is an algebraic integer for every
simple object $X$. An old conjecture says \cite[Appendix]{Kap} that the representation category
of every finite-dimensional semisimple Hopf algebra is of
Frobenius type. A classical result of Frobenius shows that if $\C$ is the category of finite-dimensional
representations of a finite group, then $\C$ is
of Frobenius type. In general, the conjecture still remains open.

In this paper we consider a class of fusion categories whose Frobenius-Perron dimension is even. The Frobenius-Perron dimensions of their simple objects are $1,2$ or $3$. We prove that such fusion category is of Frobenius type, by analyzing the structure of their Grothendieck ring.

In section \ref{sec2}, we recall the main notions and results relevant to the problem we consider. The main result is contained in Section \ref{sec3}.

\section{Preliminaries}\label{sec2}
Let $\C$ be a fusion category over $k$ and let $\Irr
(\C)$ denote the set of isomorphism classes of
simple objects of $\C$. Then $\Irr
(\C)$ is a basis of the Grothendieck ring $K_0(\C)$ of
$\C$. We use $\FPdim x$ to denote the Frobenius-Perron dimension of $x\in \Irr
(\C)$. It is the Frobenius-Perron eigenvalue of the matrix of left multiplication by $x$ in $K_0(\C)$. This extends to a ring homomorphism
$\FPdim : K_0(\C) \to \mathbb{R}$. This is the unique ring homomorphism  that takes positive values in all elements of $\Irr (\C)$. The
Frobenius-Perron dimension of $\C$ is defined by
$$\FPdim \C =\sum_{x\in \Irr(\C)} (\FPdim x)^2.$$
 A fusion category $\C$ is called integral if $\FPdim X\in \mathbb{Z}$ for all objects of $\C$. Every integral fusion category is isomorphic to the category of representations of some finite dimensional semisimple quasi-Hopf algebra \cite[Theorem 8.33]{ENO}.

Let $\C$ be a fusion category over $k$. Let $G(\C)$ denote the set of isomorphism classes of invertible objects of $\C$. Then $G(\C)$ is a subgroup of the group of units of $K_0(\C)$.

Let $X$ be an object of $\C$. Then $\FPdim X$ is defined to be the Frobenius-Perron dimension of the class of $X$ in $K_0(\C)$. As shown in \cite{ENO}, $\FPdim X \geq
1$, for all objects $X$ of $\C$. In particular, $\FPdim X = 1$ if and only
if $X$ is an invertible object.

For every $y \in K_0(\C)$, we may write $y = \sum_{x \in \Irr
(\C)} m(x, y) x$, where $m(x, y) \in \mathbb Z$. The  integer
$m(x, y)$ is called the multiplicity of $x$ in $y$. This
extends to a bilinear form $m: K_0(\C) \times K_0(\C) \to \mathbb
Z$. We then have $m(x, y)=\dim \Hom_\C(X, Y)$, where $x$ and $y$ denote the class of the objects $X$ and $Y$
of $\C$, respectively.

The following two lemmas  are restatements of  \cite[Theorem 9 and 10]{NR} in the context of fusion categories.

\medbreak\begin{lem}\label{lem1}
Let $x, y, z \in K_0(\C)$. The dual of $x$ is denoted by $x^*$. Then we have

(1)\quad $m(x,y)=m(x^*,y^*)$.

(2)\quad $m(x,yz)=m(y^*,zx^*)=m(y,xz^*)$.
\end{lem}

\medbreak\begin{lem}\label{lem2}
Let $x, y \in \Irr(\C)$. Then we have

(1)\quad For each $g \in G(\C)$, $m(g,xy)=1$ if and only if $y=x^*g$ and $0$ otherwise. In
particular, $m(g,xy)=0$ if $\FPdim x\neq \FPdim y$.

(2)\quad For all $g \in G(\C)$, $m(g,xx^{*})>0$ if and
only if $m(g,xx^{*})= 1$ if and only if $gx=x$. In particular, $G[x]=\{g\in G(\C):gx=x\}$ is a subgroup of $G(\C)$ of order at most $(\FPdim x)^2$.
\end{lem}

\medbreak Using the notations above, we have an equation
$$xx^*=\sum_{g\in G[x]}g+\sum_{y\in \Irr(\C),\FPdim y >1} m(y, xx^*) y,$$

where $x\in \Irr(\C)$.

We use $\Irr_{\alpha}(\C)$ to denote the set of
isomorphism classes of simple objects of $\C$ of Frobenius-Perron
dimension $\alpha$, where $\alpha\in \mathbb R_+$,

Let $\C$ be a fusion category. Then $\mathcal{D}\subseteq \C$ is a fusion subcategory if $\mathcal{D}$ is a full tensor
subcategory such that if $X\in \C$ is  isomorphic to a direct summand of an object of $\mathcal{D}$, then $X\in\mathcal{D}$.
If $\mathcal{D}$ is a fusion subcategory of $\C$, then $\FPdim \mathcal{D}$ divides $\FPdim \C$ \cite{ENO}. Recall that  a fusion category is called simple if it has no nontrivial proper fusion subcategories.

Fusion subcategories of $\C$ correspond to fusion
subrings of the Grothendieck ring of $K_0(\C)$, where fusion subrings means  a subring which is spanned by a subset of
$\Irr(\C)$. A subset $S\subseteq\Irr(\C)$ spans a fusion subring of
$K_0(\C)$ if and only if the product of elements of $S$ decomposes
as a sum of elements of $S$.

There is a unique largest pointed fusion subcategory of $\C$ which is generated by the group $G(\C)$ of invertible objects of $\C$. We denote this pointed fusion category by $\C_{pt}$. Moreover, the order of $G(\C)$ equals $\FPdim \C_{pt}$ and
so it divides $\FPdim \C$.

In the case when $\C$ is the representation category of a finite-dimensional semisimple Hopf algebra, the proof of the following lemma is given in \cite{NR}, while for the general case its proof is given in \cite{DNV}.

\medbreak\begin{lem}\label{lem3}
Let $x\in \Irr(\C)$. Then the following hold:

(i)\,The order of $G[x]$ divides $(\FPdim x)^2$.

(ii)\,The order of $G(\C)$ divides $n(\FPdim x)^2$, where $n$
is the number of non-isomorphic simple objects of
Frobenius-Perron dimension $\FPdim x$.
\end{lem}
\medbreak 

We call that $\C$ is of type $(d_0,n_0;
d_1,n_1;\cdots;d_s,n_s)$, where  $1 = d_0 < d_1< \cdots < d_s$ are  positive
real numbers and $n_1,n_2,\cdots,n_s$ are positive integers, if $\C$ has $n_0$ non-isomorphic simple objects of
Frobenius-Perron dimension $d_0$, $n_1$ non-isomorphic simple objects of
Frobenius-Perron dimension $d_1$, an so on.

Therefore, if $\C$ is of type $(d_0,n_0; d_1,n_1;\cdots;d_s,n_s)$,
then $n_0=|G(\C)|$ and we have an equation
\begin{equation}\label{gral} \FPdim \C = n_0 + d_1^2n_1 + \dots +
d_s^2n_s. \end{equation}

\medbreak\begin{lem}\label{lem4}
Let $\C$ be a fusion category. If there exists $x_2\in \Irr_2(\C)$ such that $x_2x_2^{*}=\1+g+x_2$ then $\C$ has a fusion subcategory of type $(1,2;2,1)$, where $g\in G(\C)$.
\end{lem}
\medbreak\begin{proof} Under the assumption above, $\{1,g,x_2\}$ spans a fusion subring of $K_0(\C)$, which corresponds to a fusion subcategory of the given type.
\end{proof}

\section{Main results}\label{sec3}

\medbreak\begin{prop}\label{pro1}
Let $\C$ be an integral fusion category such that the Frobenius-Perron dimension of every simple object is at most $3$. Suppose that there exists $x_2\in \Irr_2(\C)$ and $x_3\in \Irr_3(\C)$ such that $x_2x_2^*=\1+x_3$. Then $\C$ has a fusion subcategory of the type $(1,3;3,1)$.
\end{prop}
\medbreak\begin{proof} From $m(x_3,x_2x_2^*)=m(x_2,x_3x_2)=1$ we have $x_3x_2=x_2+w$, where $w$ is a sum of two elements of $\Irr_2(\C)$ and $m(x_2,w)=0$. Suppose that $u\in \Irr_2(\C)$ such that $m(u,w)>0$. Then $m(u,x_3x_2)=m(x_3,ux_2^*)\leq1$. Hence, $w$ is a sum of two distinct elements of $\Irr_2(\C)$. So we may write $x_3x_2=x_2+x_2'+x_2''$, where $x_2',x_2''\in \Irr_2(\C)$ are different from each other.

From $m(x_2',x_3x_2)=m(x_3,x_2'x_2^*)=1$ we have $x_2'x_2^*=g_1+x_3$ for some $1\neq g_1\in G(\C)$. Then we have $m(g_1,x_2'x_2^*)=m(x_2',g_1x_2)=1$ which means that $g_1x_2=x_2'$. So $x_2'x_2^*=g_1x_2x_2^*=g_1(\1+x_3)=g_1+x_3$ which means that $g_1x_3=x_3$. Similarly, by replacing $x_2'$ with $x_2''$, we have $g_2x_3=x_3$ for some $1\neq g_2\in G(\C)$.

We claim that $G[x_3]=\{1,g_1,g_2\}$ and hence $G[x_3]=\{1,g_1,g_1^2\}$ since $G[x_3]$ is a group.

First, $g_1\neq g_2$ otherwise $x_2'=x_2''$, a contradiction. Second, if there exists $h\in G(\C)$ such that $h\in G[x_3]$ then $hx_3=x_3$ and hence $hx_3x_2=hx_2+hx_2'+hx_2''=x_2+x_2'+x_2''$. Without loss of generality, we may assume that $hx_2=x_2'$, since we know that $hx_2\neq x_2$. Hence, $hx_2=g_1x_2$ which means that $h^{-1}g_1x_2=x_2$, so $h=g_1$ since $G[x_2]=\{\1\}$. This proves the claim by Lemma \ref{lem3}.

Now we have $x_3x_2=x_2+g_1x_2+g_1^2x_2$. Multiplying on the right by $x_2^*$, we have $x_3^2=\1+g_1+g_1^2+2x_3$. From the fact $x_3=x_3^*$ we have $(x_3g_1)^*=g_1^{-1}x_3^*=g_1^{-1}x_3=g_1^2x_3=x_3$ which means $x_3g_1=x_3$, by applying $*$ to $g_1^{-1}x_3=x_3$. Therefore, $\{1,g_1,g_1^2,x_3\}$ spans a fusion subring of $K_0(\C)$, which corresponds to a fusion subcategory of type $(1,3;3,1)$.
\end{proof}

\begin{rem}
In the semisimple Hopf algebra setting, a similar result was proved as a part of \cite[Theorem 11]{NR}. The two proofs are slightly different.
\end{rem}

\medbreak\begin{cor}\label{cor1}
Let $\C$ be an integral fusion category such that the Frobenius-Perron dimension of every simple object is at most $3$. Suppose that  $\Irr_2(\C)\neq \varnothing $. Then at least one of the followings holds:

(1)\, $\C$ has a fusion subcategory of type $(1,3;3,1)$.

(2)\, There exists $g\in G(\C)$ of order $2$ and $\C$ has a fusion subcategory of type $(1,n_0;2,n_1)$, where $n_0=|G(\C)|$ and $n_1=|\Irr_2(\C)|$.

In particular, $2$ divides the Frobenius-Perron dimension of $\C$.
\end{cor}

\medbreak
 \begin{proof}If there exists $x_2\in \Irr_2(\C)$ such that $x_2x_2^*=\1+x_3$ for some $x_3\in\Irr_3(\C)$, then part (1) follows from Proposition \ref{pro1}, otherwise part (2) follows from \cite[Lemma 3.2(a)]{DNV}.
\end{proof}

\medbreak\begin{prop}\label{pro2} Let $\C$ be an integral fusion category such that the Frobenius-Perron dimension of every simple object is at most $3$. Suppose that $\Irr_3(\C)\neq \varnothing $. Then at least one of the followings holds:

(1)\,$\C$ has a pointed fusion subcategory of Frobenius-Perron dimension $3$.

(2)\,$\C$ has a fusion subcategory of type $(1,3;3,1)$, $(1,2;2,1)$ or $(1,2;2,4)$.

In particular, $3$ divides the Frobenius-Perron dimension of $\C$.
\end{prop}
\medbreak\begin{proof} Let $x_3\in \Irr_3(\C)$. By Lemma \ref{lem3}, the order of $G[x_3]$ is $1,3$ or $9$. If $|G[x_3]|=3$ or $9$ then part (1) holds true. So it remains to consider the case when $|G[x_3]|=1$.

There are two possible decompositions of $x_3x_3^*$:
$$x_3x_3^*=\1+x_2+x_3'+x_3''\mbox{\quad or \quad}x_3x_3^*=\1+a_1+a_2+a_3+a_4,$$
where $x_2,a_i(i=1,2,3,4)\in\Irr_2(\C)$ and $x_3',x_3''\in\Irr_3(\C)$. In particular, $x_2$ is self-dual since $x_3x_3^*$ is self-dual.

Suppose that the first decomposition of $x_3x_3^*$ holds true. From $1=m(x_2,x_3x_3^*)=m(x_3,x_2x_3)$ we have $x_2x_3=x_3+u$, where $\FPdim u=3$. Since $\FPdim x_2\neq \FPdim x_3$, $u$ can not contain elements from $G(\C)$. Therefore, $u\in\Irr_3(\C)$ and $u\neq x_3$. From $1=m(u,x_2x_3)=m(x_2,ux_3^*)$ we have $ux_3^*=x_2+w$, where $\FPdim w=7$. Multiplying $x_2x_3=x_3+u$ on the right by $x_3^*$, we have
\begin{equation}\label{equ1}x_2^2+x_2x_3'+x_2x_3''=\1+x_2+x_3'+x_3''+w.\end{equation}

We shall consider the three possible decompositions of $x_2^2=x_2x_2^*$.

Assume that $x_2^2=\1+v$, where $v\in \Irr_3(\C)$ . Then Proposition \ref{pro1} shows that $\C$ has a fusion subcategory of type $(1,3;3,1)$.

Assume that $x_2^2=\1+g+v$, where $g\in G(\C)$ and $v\in \Irr_2(\C)$. If $v\neq x_2$ then equation (\ref{equ1}) shows that $x_2$ lies in $x_2x_3'$ or $x_2x_3''$. Then $m(x_2,x_2y)\geq 1$, where $y$ is $x_3'$ or $x_3''$. This implies that $m(x_2,y^*x_2)=m(y^*,x_2^2)\geq 1$, which contradicts with the assumption  on the decompositions of $x_2^2$. Therefore, $x_2^2=\1+g+x_2$ and hence $\C$ has a fusion subcategory of type $(1,2;2,1)$ by Lemma \ref{lem4}.

Assume that $x_2^2=\1+g_1+g_2+g_3$, where $g_i\in G(\C),i=1,2,3$. In this case, $x_2$ must lie in the decomposition of $x_2x_3'$ or $x_2x_3''$ by equation (\ref{equ1}). A similar argument as in the paragraph above shows that it is impossible. So we get a contradiction in this case.

Suppose that the second decomposition of $x_3x_3^*$ holds true. From $m(a_i,x_3x_3^*)=m(x_3,a_ix_3)\geq 1$ we have $a_ix_3=x_3+w_i$, where $w_i\in \Irr_3(\C)$ and $i=1,2,3$ or $4$.

Assume that $w_i=x_3$. Then $a_ix_3=2x_3$ and hence $x_3x_3^*=\1+2a_i+b'+c'$, where $b',c'\in \Irr_2(\C)$ are different from $a_i$. Multiplying $a_ix_3=2x_3$ on the right by $x_3^*$, we have
\begin{equation}\label{equ2}2a_i^2+a_ib'+a_ic'=2\cdot\1+3a_i+2b'+2c'.\end{equation}
From the equality above we first know that $b'\neq c'$, since the right hand side is not divisible by 2, hence $\1$ must lie in the decomposition of $a_i^2$. It follows that $a_i$ is self-dual and $G[a_i]=\{\1\}$. We then may write $a_i^2=\1+w$, where $w\in \Irr_3(\C)$. But it is impossible because the right side of equality $(\ref{equ2})$ does not contain elements of $\Irr_3(\C)$. Therefore, $w\neq x_3$ and hence $m(a_i,x_3x_3^*)=1$. Hence, the multiplicity of $a_1,a_2,a_3$ and $a_4$ in $x_3x_3^*$ is $1$, respectively. In other words, $a_1,a_2,a_3,a_4$ are distinct.

\medbreak By Proposition \ref{pro1}, if there exists $i$ such that $a_ia_i^*=\1+v$ for some $v\in\Irr_3(\C)$ then $\C$ has a fusion subcategory of type $(1,3;3,1)$.  We are done in this case. We therefore assume that the order of $\Pic[a_i]$ is greater than  $1$ for all $i$. In addition, Lemma \ref{lem3} shows that in this case the order of $\Pic[a_i]$ is $2$ or $4$ for all $i$.

Let $1\neq g,h\in \Pic[a_i]$. Then $ga_i=ha_i=a_i$. So
\begin{eqnarray*}
&&a_ix_3=ga_ix_3=g(x_3+w_i)=gx_3+gw_i=x_3+w_i, \mbox{and}\\
&&a_ix_3=ha_ix_3=h(x_3+w_i)=hx_3+hw_i=x_3+w_i.
 \end{eqnarray*}
Hence $gx_3=w_i=hx_3$. Since $\Pic[x_3]=\{\1\}$, we have $g=h$. Therefore,  the order of $\Pic[a_i]$ is $2$ for all $i$.

\medbreak We may write $a_ia_i^*=\1+g_i+b_i$ for some $g_i\in \Pic(\C)$ and $b_i\in \Irr_2(\C)$. Since $a_i^*$ also lies in the decomposition of $x_3x_3^*$, we have $m(a_i^*,x_3x_3^*)=m(x_3,a_i^*x_3)=1$, which implies that $a_i^*x_3=x_3+w_i'$, where $w_i'\in \Irr_3(\C)$. In addition, from $m(w_i',a_i^*x_3)=m((w_i')^*,x_3^*a_i)=m(x_3^*,(w_i')^*a_i^*)=m(x_3,a_iw_i')=1$,
we have $a_iw_i'=x_3+w_i''$, where $w_i''\in \Irr_3(\C)$. Multiplying $a_i^*x_3=x_3+w_i'$ on the left by $a_i$, we have
$$(\1+g_i+b_i)x_3=x_3+g_ix_3+b_ix_3=2x_3+w_i+w_i'',$$
which implies that $g_ix_3+b_ix_3=x_3+w_i+w_i''$. This means that $m(x_3,b_ix_3)=m(b_i,x_3x_3^*)=1$. Hence, $b_i$ also lies in the decomposition of $x_3x_3^*$ for all $i$.

We claim that there exists $g\in \Pic(\C)$ such that $\Pic[a_i]=\{\1,g\}$ for all $1\leq i\leq 4$. In fact, in our situation, $a_1a_i^*$ is not irreducible for $i=2,3,4$. Then \cite[Lemma 2.5]{DNV} shows that $a_1a_1^*$ and $a_ia_i^*$ must contain a common nontrivial irreducible component. If $\Pic[a_1]\neq \Pic[a_i]$ then we can write
$$a_1a_1^*=\1+g+b\quad\mbox{and} \quad a_ia_i^*=\1+h+b,$$
where $g,h\in\Pic(\C)$ and $b\in\Irr_2(\C)$. Then the first one implies that $\{1,g\}\subseteq G[b]$, and the latter one implies that $\{1,h\}\subseteq G[b]$. Furthermore, the discussion in the paragraph above shows that $b\in \{a_1,a_2,a_3,a_4\}$, and hence the order of $\Pic[b]$ is $2$. Therefore, we get the contradiction that $\{1,g\}=\{1,h\}$.

As a conclusion, we obtain that $ga_i=a_i$ for all $i$. Since $\{a_1,a_2,a_3,a_4\}=\{a_1^*,a_2^*,a_3^*,a_4^*\}$, we also obtain that $ga_i^*=a_i^*$. Hence, $a_ig=a_i$ for all $i$.

From $a_ix_3=x_3+w_i$ and $ga_i=a_i$, we obtain that $w_i=gx_3$. Multiplying $a_ix_3=x_3+gx_3$ on the right by $x_3^*$, we have
\begin{eqnarray*}
a_i(\1+a_1+a_2+a_3+a_4)&=&a_i+a_ia_1+a_ia_2+a_ia_3+a_ia_4\\
&=&\1+a_1+a_2+a_3+a_4+g(\1+a_1+a_2+a_3+a_4)\\
&=&\1+g+2a_1+2a_2+2a_3+2a_4.
\end{eqnarray*}
Therefore, $a_ia_1, a_ia_2, a_ia_3$ and $a_ia_4$ are the sums of elements of $\{\1,g,a_1,a_2,a_3,a_4\}$, respectively. Combining this result with the fact that $ga_i=a_ig=a_i$, we obtain that $\{\1,g,a_1,a_2,a_3,a_4\}$ spans a fusion subring of $K_0(\C)$. It follows that $\C$ has a fusion subcategory of type $(1,2;2,4)$. This completes the proof.
\end{proof}
 Combining Corollary \ref{cor1} with Proposition \ref{pro2}, we obtain our main result.
\begin{thm}
Let $\C$ be an integral fusion category such that the Frobenius-Perron dimension of every simple object is at most $3$.  Then $\C$ is of Frobenius type and $\C$ is not simple.
\end{thm}

\begin{rem}
In the semisimple Hopf algebra setting, there are two relevant results: Let $H$ be a finite-dimensional semisimple Hopf algebra over $k$. In \cite{NR}, the authors proved that if $H$ has a simple module of dimension $2$ then $2$ divides the dimension of $H$. In \cite{Bur,Don,KSZ}, the authors proved that if dim$H$ is odd and $H$ has a simple module of dimension $3$ then $3$ divides the dimension of $H$. However, it is not known whether $3$ divides dim$H$ when dim$H$ is even and $H$ has a simple module of dimension $3$. Our present work gives a partial answer to this question.
\end{rem}

\end{document}